\documentstyle[12pt]{article}
\newtheorem{theorem}{Theorem}
\newtheorem{lemma}{Lemma}
\newtheorem{proposition}{Proposition}
\newtheorem{remark}{Remark}
\newtheorem{corollary}{Corollary}
\newtheorem{definition}{Definition}
\newtheorem{example}{Example}

\newlength{\rig}
\newlength{\rigg}
\newlength{\hei}

\newcommand{\func}[1]{{\rm #1} \,}

\newcommand{\limfunc}[1]{{\rm{#1}}}
\newcommand{\text}[1]{{\mbox {#1}}}
\begin{document}

\setcounter{page}{1}

\author{{\bf Anton Savin\footnotemark\rule{5pt}{0pt} and
\addtocounter{footnote}{-1}
Boris Sternin\footnotemark}%
\addtocounter{footnote}{-1}\thanks{
Supported by the Russian Foundation for Basic Research under
grants N 97-01-00703 and 97-02-16722a and by Arbeitsgruppe Partielle
Differentialgleichungen und Komplexe Analysis, Institut f\"ur Mathematik,
Universit\"at Potsdam, and Soros Foundation.}
\\[3mm]
Moscow State University\\[3mm]
e-mail: antonsavin@mtu-net.ru\\[3mm]
e-mail: sternin@math.uni-potsdam.de}
\title{\bf Elliptic Operators in Odd Subspaces}
\date{}
\maketitle

\vspace{15mm}

\begin{abstract}
An elliptic theory is constructed for operators acting in subspaces
defined via odd pseudodifferential projections. Index formulas are
obtained for operators on compact manifolds without boundary and for
general boundary value problems. A connection with Gilkey's theory of
$\eta$-invariants is established.
\end{abstract}

\vspace{15mm}


\vspace{15mm}

{\bf 1991 AMS classification}: Primary 58G03, Secondary 58G10, 58G12,
58G25, 19K56

\vfill

\newpage

\tableofcontents

\vspace{0.3cm}

\section*{Introduction}
\addcontentsline{toc}{section}{Introduction}

The present paper is a continuation of \cite{SaSt1},
where we studied subspaces
defined by pseudodifferential projections
and constructed a homotopy invariant functional
\begin{equation}
\label{al1}
d: {\rm \widehat{Even}}(M^{odd})\longrightarrow{\bf Z}\left[\frac 12\right]
\end{equation}
on the class of even subspaces on odd-dimensional manifolds. The
functional (\ref{al1}) is a generalization of the notion of dimension of a
finite-dimensional vector space. In terms of this functional, we obtained
index formulas for general boundary value problems with boundary
conditions ranging in even subspaces and for elliptic operators in such
subspaces.

Using Gilkey's results \cite{Gil7}, we found that (under certain
unessential restrictions) the functional $d$ coincides with the
$\eta$-invariant of Atiyah--Patodi--Singer for some self-adjoint elliptic
operator $A$ of even order. Let us point out that in this case the
$\eta$-invariant is invariant with respect to homotopies of self-adjoint
operators modulo integer jumps related to the spectral flow (see
\cite{APS2}).

On the other hand, it can be shown that the $\eta$-invariant is also a
homotopy invariant (modulo integer jumps) for the class of odd-order
operators on even-dimensional manifolds (see \cite{Gil8}). Operators of
this kind and their $\eta$-invariants found numerous applications in
differential geometry and topology of manifolds (e.g., see
\cite{BaGi1,Gil10,Stol1,LiZh1,Liu1}).

Thus, there arises a natural question of defining the functional $d$ on
the subspaces corresponding to elliptic self-adjoint operators of odd
order.

The present paper deals with the construction of the corresponding
functional on the set of ``odd'' subspaces, i.e., subspaces corresponding
to self-adjoint operators of odd order. As in \cite{SaSt1}, in terms of
the functional $d$ we obtain index formulas for general boundary value
problems in odd subspaces, as well as an index formula for elliptic
operators acting in subspaces.

We also prove the equality of the functional $d$ and the $\eta$-invariant.
This equality, in particular, proves Gilkey's conjecture \cite{Gil8} that
the $\eta$-invariant is always a dyadic rational.

The subject of this paper is largely analogous to the subject of our previous
paper \cite{SaSt1}. However, a number of constructions here are totally
different. Because of this we are rather brief in the part of the paper
that parallels \cite{SaSt1}, concentrating on the problems radically
different from those considered previously.

Let us briefly describe the contents of the paper.

In the first section the notion of an odd subspace is introduced. Here
we also show the relation of such subspaces to the spectral subspaces of
elliptic self-adjoint operators with odd symbols. It is shown that an
arbitrary odd subspace can be transformed into a spectral subspace of
this form. In the second section we define the functional $d$ on odd
subspaces on
even-dimensional manifolds. In terms of this functional, in sections 3 and
5 homotopy invariant index formulas are obtained for elliptic operators
acting in the subspaces and for the general boundary value problems. In
the fourth section, the equality of the functional $d$ and the
$\eta$-invariant of Atiyah-Patodi-Singer is established. In the next
section, an example is considered in which both sides of the index formula
for an operator in subspaces are computed independently and their equality
is verified. At the end of the paper we place a short appendix explaining
the notion of relative index of projections and of subspaces defined by
projections.

\noindent
{\bf Preliminary publications.}
The results of the paper were reported at the international conference
``Operator algebras and asymptotics on manifolds with singularities,''
Warsaw, Poland \cite{SaSt3}.
The paper appeared as a preprint in \cite{SaSt2}

\section{Odd vector bundles. Odd subspaces}

Let $M$ be a closed smooth manifold and $E\in{\rm Vect}\left( M\right) $
a vector bundle. Consider the involution
\[
\alpha :T^{*}M\rightarrow T^{*}M,\quad \alpha \left(x, \xi \right) =(x,-\xi)
\]
on the cotangent bundle.

\begin{definition}
{\em A subbundle} $L\subset \pi ^{*}E$, $\pi :S^{*}M\rightarrow M,$ {\em
is said to be} odd {\em with respect to the involution $\alpha$ if it is a
complement of $\alpha^{*}L$ in $\pi^{*}E$, that is,}
\begin{equation}
L\oplus \alpha ^{*}L=\pi ^{*}E.\quad   \label{comple}
\end{equation}
\end{definition}

Let $L\subset \pi ^{*}E$ be an odd bundle. The
{\em odd projection\/} on $L$ in $\pi ^{*}E$ is the projection
\begin{equation}
\label{twostars}
p_L:\pi ^{*}E\rightarrow \pi ^{*}E,
\end{equation}
on $L$ along $\alpha ^{*}L.$ The condition that a projection $p_L$ is odd
can be represented in the form
\begin{equation}
\label{sstar}
p_L+\alpha ^{*}p_L=Id.
\end{equation}

\begin{remark}
{\em The minimum possible dimension of an odd subbundle $L\subset \pi ^{*}
E$ is very large as compared with the dimension of $M$. More precisely, on
a manifold $M$ of dimension $n$ one has the following relation:}
\begin{equation}
\label{onestar}
\left.
\begin{array}{l}
n-1=2k \\
n-1=2k+1
\end{array}
\right\} \Rightarrow \dim L\,\text{\em is a multiple of \ }2^{k-1}.
\end{equation}
\end{remark}

Indeed, for an odd projection $p_L$ the matrix-valued function
\[
2p_L-Id
\]
is odd one on the sphere $S^{n-1}$.
Now relation (\ref{onestar}) follows readily
from Theorem 1 in \cite{Gil9}, which states that an odd function
on the sphere ranging in invertible matrices obeys (\ref{onestar}).

\begin{example}
\label{thefirst}
{\em
Let $M=S^1$, and let $E$ be the trivial line bundle,
$E={\bf C}$. The cotangent sphere bundle in this case is the disjoint
union of two circles,
$$
S^*S^1=S^1_+\bigsqcup S^1_-,\quad S^1_\pm=
\left\{
(x,\xi)\in T^*S^1\;|\; \xi=\pm 1
\right\}.
$$
Consider the subbundle $L\subset \pi ^{*}E,\;\pi :S^{*}S^1\to S^1$,
such that $L={\bf C}$ over the subspace $S^1_+$ and
$L=\{0\}$ over $S^1_-$. This bundle is obviously odd.
}
\end{example}

If $E$ is equipped with a metric, then
it is useful to consider the class of bundles $L$ for which
the direct sum (\ref{comple}) is orthogonal.
Bundles $L$ of this type are said to be {\em orthogonally odd\/}.

\begin{remark}
{\em If
\[
\sigma :\pi ^{*}E\longrightarrow \pi ^{*}F
\]
is an isomorphism of vector bundles and $\sigma$ is even in the sense that
$$
\alpha ^{*}\sigma =\sigma,
$$
then each odd subbundle $L\subset \pi ^{*}E$ is taken by $\sigma$
into an odd subbundle $\sigma L\subset \pi ^{*}F$. }
\end{remark}

\begin{definition}
{\em A subspace $\widehat{L}\subset C^\infty \left(
M,E\right) $ is said to be} admissible\/ {\em if
\[
\widehat{L}=\func{Im}P
\]
for some pseudodifferential projection $P$ of order zero that acts in the
space of sections of the vector bundle $E$. }
\end{definition}

The subbundle $L=\func{Im}\sigma \left( P\right) $ is called the
(principal) {\em symbol of the subspace\/} $\widehat{L}$ (it is
independent of the choice of projection on $\widehat{L}$). A subspace
$\widehat{L}$ is called {\em odd\/} if its symbol $L$ is odd. The
semigroup of odd subspaces (with respect to the direct sum of the ambient
spaces) is denoted by $\widehat{\limfunc{Odd}}\left(M\right)$.

Odd subspaces can be constructed by means of odd projections
(i.e. pseudodifferential projections with principal symbols
satisfying (\ref{sstar})) or with the help of self-adjoint
elliptic operators, by virtue of the following proposition.

\begin{proposition}
\label{opisub}
Let $A$ be an elliptic self-adjoint operator on $M$ with odd principal
symbol:
\begin{equation}
\label{anti}
\alpha^{*}\sigma(A)=-\sigma(A).
\end{equation}
 Then the subspace $\widehat{L}_+(A)$ generated by the
eigenvectors of $A$ corresponding to nonnegative eigenvalues is odd.
\end{proposition}
\begin{remark}
{\em
The condition mentioned in Proposition \ref{opisub}
is satisfied, for example, for all differential operators of odd order.
}
\end{remark}

{\em Proof\/}. It is well known (e.g., see \cite{APS1}) that the principal
symbol of the orthogonal projection on the nonnegative spectral subspace $
\widehat{L}_+(A)$ is equal to the projection on the nonnegative spectral
subspace ${L}_+(A)$ of the principal symbol $\sigma(A)$ along the
corresponding negative subspace. By virtue of condition (\ref{anti}), the
subbundle ${L}_+(A) \subset\pi^{*}E$ is odd. This proves the proposition.

\begin{example}
{\em
Consider the Hardy space, i.e., the subspace
$$
\widehat{L}\subset C^\infty(S^1)
$$
of boundary values of holomorphic functions in the unit disk
$$
D\subset {\bf C},\quad \partial D = S^1.
$$
A simple calculation with Fourier series shows that the Hardy space
$\widehat{L}$ coincides with the nonnegative spectral subspace of an
elliptic (essentially) self-adjoint operator
$A=-i\frac{d}{d\varphi}$ on the circle $S^1$. In particular, the symbol
of the Hardy space is equal to
$$
L_{\varphi,\xi}=
\left\{
\begin{array}{ll}
{\bf C} & \mbox{ for }\quad \xi=1, \\
0   &\mbox{ for } \quad \xi=-1.
\end{array}
\right.
$$
It is clear that the symbol of the Hardy space is just the
odd bundle from Example \ref{thefirst}.
}
\end{example}

General  odd subspaces can be reduced to the orthogonally odd subspaces
in the sense of the following proposition.

\begin{proposition}
\label{reduc}For an odd subspace $\widehat{L}$, there exists
an invertible pseudodifferential operator $U$
\[
U:C^\infty \left( M,E\right) \rightarrow C^\infty \left( M,E\right)
\]
with even symbol such that the symbol of the subspace $U\widehat{L}$ is
orthogonally odd.
\end{proposition}

{\em Proof\/}. Let us consider two metrics in the bundle
$\pi ^{*}E$, a metric $g_0$ induced by an arbitrary metric on $E$
and an arbitrary metric $g_1$
even with respect to the involution
$\alpha$ and such that the sum ${L}\oplus \alpha ^{*}{L}$
is an orthogonal sum. For the linear homotopy
\[
g_t=\left( 1-t\right) g_0+tg_1,
\]
of these metrics, we take a covering homotopy of unitary transformations
\[
u_t:\left( \pi ^{*}E,g_0\right) \rightarrow \left( \pi ^{*}E,g_t\right)
\]
and consider an operator $U_1$ with principal symbol $u_1.$
It can be chosen to be invertible, since its symbol
is homotopic to unity. The inverse $U=\left( U_1\right) ^{-1}$
is the desired operator.

Proposition \ref{reduc} permits us to extend results valid for orthogonally
odd subspaces to the general case.

\begin{proposition}
An odd subspace $\widehat{L}$ can be specified by an odd projection
$P_{\widehat{L}}$. In particular, an odd subspace has an odd complement.
\end{proposition}

{\em Proof\/}. For the case in which $\widehat{L}$ is orthogonally odd, we
can take $P_{\widehat{L}}$ to be the orthogonal projection operator on
$\widehat{L}$. The general case follows by Proposition \ref{reduc}.

\section{The structure of odd subspaces. The dimension functional
\label{dime}}

In this section we obtain a stable classification of odd subspaces on
even-dimensional manifolds with respect to invertible transformations
given by even elliptic operators, i.e., pseudodifferential operators with
even principal symbols. On the basis of this classification, we define an
analog of the notion of dimension for odd subspaces in the end of the
section.

\begin{definition}
{\em Two subspaces $\widehat{L}_1\subset C^\infty \left(
M,E\right) $ and $\widehat{L}_2\subset C^\infty \left( M,F\right) $
are called} equivalent\/ {\em if there exists an even invertible operator
\[
U:C^\infty \left( M,E\right) \rightarrow C^\infty \left( M,F\right)
\]
transforming $\widehat{L}_1$ into $\widehat{L}_2$. }
\end{definition}

The abelian semigroup of equivalence classes of odd subspaces is denoted
by $\left. \widehat{\rm Odd}\right/ \sim $. Let us denote the
corresponding Grothendieck group by $K\left(\left. \widehat{\rm Odd}\right
/ \sim \right)$. We also consider the following equivalence relation in
the usual semigroup ${\rm Vect}\left( M\right)$ of vector bundles:
\[
E\sim F
\]
if there exists an even elliptic symbol
\[
\sigma :\pi ^{*}E\rightarrow \pi ^{*}F
\]
with trivial (topological) index. The corresponding Grothendieck group
is denoted by $K({\rm Vect}(M)/\sim)$.\footnote{This equivalence
relation induces also an equivalence on the Grothendieck
group $K\left( M\right) .$ There is, in particular, a natural
isomorphism of the groups
\[
K\left( \left. {\rm Vect}\left( M\right) \right/ \sim \right)
\quad \text{and }
\left. K\left( M\right) \right/ \sim .
\]
}

Consider the sequence
\begin{equation}
0\rightarrow {\bf{}Z}\stackrel{i}{\rightarrow }
K\left( \left. \widehat{\rm Odd}%
\right/ \sim \right) \stackrel{j}{\rightarrow } K\left(
{\rm Vect}(M)/\sim\right)
\rightarrow 0.  \label{seq1}
\end{equation}
Here $i$ takes each integer $k$ to the element
\begin{equation}
i\left( k\right) =\left[ \widehat{L}+k\right] -\left[ \widehat{L}\right] ,
\label{easy}
\end{equation}
where $\widehat{L}\subset C^\infty \left( M,E\right) $ is an
odd subspace and $\widehat{L}+k\subset C^\infty \left(
M,E\right) $ is the sum of $\widehat{L}$ and an arbitrary subspace
of dimension $k$ in the complement of $\widehat{L}$.
We note that the difference (\ref{easy})
is independent of the choice of $\widehat{L}$.
The homomorphism $j$ is induced by the mapping that takes each
subspace $\widehat{L}\subset C^\infty \left(M,E\right)$
to the bundle $E$.

Finally, note that the sequence (\ref{seq1}) is a complex: $ji=0.$

\begin{theorem}\label{thm101}
On an even-dimensional manifold $M$, the sequence {\em (\ref{seq1})} is $
2$-exact. More precisely, the mapping $i$ is an embedding, while the
groups $\limfunc{coker}j$ and $\left. \ker j\right/ \limfunc{im}\,i$ are $
2$-primary groups, i.e. torsion groups with orders powers of 2.
\end{theorem}

\begin{remark}
{\em In terms of localization (see \cite{Sul1}, Sect. 1),
Theorem~\ref{thm101} states that the localization of the sequence outside
of two is exact, that is, the following sequence is exact
\begin{equation}
0\rightarrow {\bf{}Z}\left[ \frac 12\right] \stackrel{i}{\rightarrow }%
K\left( \left. \widehat{\rm Odd}\right/ \sim \right) \otimes {\bf{}Z}\left[
\frac 12\right] \stackrel{j}{\rightarrow }
K\left({\rm Vect}(M)/\sim\right)
\otimes {\bf{}Z}\left[ \frac 12\right] \rightarrow 0,
\label{loca}
\end{equation}
where ${\bf{}Z}\left[ \frac 12\right] $ is the group of dyadic rationals
and the tensor products are taken over ${\bf Z}$. }
\end{remark}

\noindent The proof of the theorem is divided into several
lemmas.

\begin{lemma} \label{compact}
The mapping $i$ is an isomorphism onto the subgroup $K_0\left( \left.
\widehat{\rm Odd}\right/ \sim \right) $ generated by the differences of
subspaces with equal symbols. In particular, this subgroup carries a
well-defined relative index of subspaces (see Appendix).
\end{lemma}

{\em Proof\/}. By construction, we obtain
$$
\func{Im}i\subset K_0\left( \left. \widehat{ Odd}\right/ \sim \right)
$$
(see (\ref{easy})). Let us now verify that $i$ is an embedding. Let
$i\left( k\right) =0$; then the subspaces $\widehat{L}+k$ and $\widehat{L}
$ are equivalent,\footnote{Here and in what follows, standard
considerations related to the stabilization of elements in Grothendieck
groups are omitted. Whenever additional constructions are needed, they are
discussed separatelly.} i.e., they are transformed into each other by an
even invertible operator $D$.

Consider the operator $D:\widehat{L}+k\rightarrow \widehat{L}$, which
is an analog of Toeplitz operators
(e.g., see \cite{BaDo1}), i.e., an operator, acting in the subspaces
with equal symbols.
Since $D$ is invertible, we obtain
\begin{equation}
{\rm ind}\left( D,\widehat{L}+k,\widehat{L}\right) =0.  \label{A1}
\end{equation}
At the same time, by applying the index formula for Toeplitz
operators\footnote{For an elliptic operator
\[
D:\widehat{L}_1\rightarrow \widehat{L}_2
\]
acting in subspaces with equal symbols $L_1=L_2$,
this formula states
\[
{\rm ind}\left( D,\widehat{L}_1,\widehat{L}_2\right) ={\rm ind}\left( \left. \sigma
\left( D\right) \right| _L\oplus 1_{L^{\perp }}\right) +{\rm ind}\left( \widehat{L}%
_1,\widehat{L}_2\right) .
\]
This is an easy exercise on the logarithmic property of the index.}
to the operator $D$, we obtain
\begin{equation}
{\rm ind}\left( D,\widehat{L}+k,\widehat{L}\right) ={\rm ind}\left(
\left. \sigma \left(
D\right) \right| _L\oplus 1_{\alpha ^{*}L}\right) +k.  \label{A2}
\end{equation}
Let us show that
\begin{equation}
{\rm ind}\left( \left. \sigma \left( D\right) \right| _L\oplus 1_{\alpha
^{*}L}\right) =\frac 12{\rm ind}D
\label{A3}
\end{equation}
for an even operator $D$ on an even-dimensional manifold.
Indeed, the cohomological Atiyah-Singer formula, together with
the fact that the involution $\alpha$ preserves the orientation
of the cotangent space $T^{*}M$, implies
\[
{\rm ind}\left( \left. \sigma \left( D\right) \right| _L\oplus 1_{\alpha
^{*}L}\right) ={\rm ind}\left( \alpha ^{*}\left( \left. \sigma \left( D\right)
\right| _L\oplus 1_{\alpha ^{*}L}\right) \right).
\]
On the other hand, since $D$ is even and $L$ is odd, we obtain
\[
\alpha ^{*}\left( \left. \sigma \left( D\right) \right| _L\oplus 1_{\alpha
^{*}L}\right) =\left. \sigma \left( D\right) \right| _{\alpha ^{*}L}\oplus
1_L.
\]
Hence,
\[
2\cdot{\rm ind}\left( \left. \sigma \left( D\right) \right| _L\oplus 1_{\alpha
^{*}L}\right) ={\rm ind}\left( \left. \sigma \left( D\right) \right| _{L\oplus
\alpha ^{*}L}\oplus 1_{\alpha ^{*}L\oplus L}\right) ={\rm ind}\,D.
\]

By substituting (\ref{A1}), (\ref{A3}) into (\ref{A2}), and by taking into
account the invertibility of the operator $D$, we obtain $k=0$, as
desired.

Let us check the opposite inclusion
$K_0\left( \left. \widehat{\rm Odd}\right/ \sim \right)
\subset \func{Im}i$. Consider two odd subspaces $\widehat{L}_{1,2}$
with equal symbols. The corresponding orthogonal projections
$P_{1,2}$ are
pseudodifferential operators with equal principal symbols.
Hence, it follows from the homotopy classification of projections with
equal principal symbols that
$P_1$ is homotopic to $P_2$ modulo a finite rank projection
(see \cite{Wojc1}, cf. \cite{EpMe1}). More precisely, the projections
become homotopic after we add a finite rank projection to one of them.
For a homotopy $P_t$ of projections, consider the invertible operator
$U_t$ given by the solution of the Cauchy problem
\[
\stackrel{.}{U_t}=\left[ \stackrel{.}{P_t},P_t\right] U_t,\quad U_1={\rm id}.
\]
An explicit check shows that $\func{Im}P_t=\func{Im}U_tP_1.$
This implies that the subspace $\widehat{L}_1$ is equivalent to
$\widehat{L}_2$ modulo a finite-dimensional subspace. Q.E.D.

\begin{lemma}
\label{coke}The quotient group
$\limfunc{coker} j=\left.
  K\left({\rm Vect}(M)/\sim\right)
\right/ \func{Im}j$ is a 2-primary group, i.e.,
for an arbitrary vector bundle
$E\in {\rm Vect}\left( M\right) $ and for a sufficiently large positive
integer $N$ there exists an odd subbundle
in the direct sum
\[
2^N\pi ^{*}E,\quad
\mbox{ for } 2^NE=\underbrace{E\oplus\ldots\oplus E}_{2^N
\mbox{\footnotesize \em times}}.
\]
\end{lemma}

{\em Proof\/}. Let us verify the validity of this statement for the
trivial bundle. We embed the manifold $M$ in Euclidean space of
dimension $N$. The cotangent bundle is then embedded in the trivial
bundle,
\begin{equation}
\label{norml}
T^{*}M\oplus \nu ={\bf{}R}^N,
\end{equation}
where $\nu$ is the normal bundle to $M$ in ${\bf R}^N$. Consider the
Clifford algebra ${\rm Cl}\left( {\bf{}C}^N\right) $ of the Hermitian
space ${\bf{}C}^N$ (e.g., see \cite{LaMi1}). Obviously, ${\bf{}C}^N$ (and
hence ${\bf{}R}^N$) acts on the space ${\rm Cl}\left( {\bf{}C}^N\right) $,
$\dim {\rm Cl}\left( {\bf{}C}^N\right) =2^N$ by the Clifford
multiplication. Moreover, using the embedding (\ref{norml}), we obtain the
homomorphism
$$
T^{*}M\stackrel{cl}{\rightarrow } {\rm End}\left( {\rm
Cl}\left( {\bf{}C}^N\right)\right)
$$
of vector bundles given by the Clifford multiplication. For the covariant
differential $\bigtriangledown=d\otimes 1_{2^N}$ in the vector bundle with
fiber ${\rm Cl}\left({\bf{}C}^N\right) $, consider the self-adjoint Dirac
operator (e.g., see \cite{BaDo1})
\[
C^\infty \left( M,{\rm Cl}\left( {\bf{}C}^N\right) \right) \stackrel{%
\bigtriangledown }{\rightarrow }C^\infty \left( M,T^{*}M\otimes {\rm Cl}\left(
{\bf{}C}^N\right) \right) \stackrel{cl}{\rightarrow }C^\infty \left(
M,{\rm Cl}\left( {\bf{}C}^N\right) \right) .
\]
The positive subspace of its principal symbol is odd. This proves the
lemma for the case of a trivial bundle. The proof for the general case
$E\in{\rm Vect}(M)$ is reduced to the one just described by means of the
tensor product by the bundle $E$.

\begin{lemma}
\label{fund}For an arbitrary odd bundle $L$ and
sufficiently large $N$, there exists an isomorphism
\[
2^NL\sim 2^N\alpha ^{*}L
\]
of vector bundles over the space $S^{*}M$ of cotangent spheres. Moreover,
there exists an even isomorphism with zero (topological) index
transforming the subbundle $2^NL\subset 2^N\pi^*E$
into $2^N\alpha ^{*}L$.
\end{lemma}

{\em Proof\/}. It follows from the stabilization theorem \footnote{This
theorem claims that the equality $\left[ E_1\right] =\left[ E_2\right] \in
K\left( X\right) $ of two elements of the $K$-group implies the
isomorphism of vector bundles $E_1$ and $E_2$ provided the dimension of
these bundles is large enough compared with the dimension of $X$.} (e.g.,
see \cite{Hus1}) that to prove the first part of the proposition, it
suffices to check the equality of the elements $2^NL$ and $2^N\alpha
^{*}L$ in the group $K\left(S^{*}M\right) $.

Consider the natural projection $p:S^{*}M\longrightarrow P^{*}M$ of the
bundle of cotangent spheres on its projectivization. Over a point $x$ of
the manifold $M$, the corresponding standard mapping
\[
p_x:S^{2n-1}\longrightarrow {\bf{}RP}^{2n-1}
\]
induces a 2-isomorphism in $K$-theory
\[
p^{*}:K^{*}\left( {\bf{}RP}^{2n-1}\right) \longrightarrow K^{*}\left(
S^{2n-1}\right) ,\quad \ker p^{*}={\bf{}Z}_{2^{n-1}},\limfunc{coker}p^{*}=%
{\bf{}Z}_2.
\]
Now by the Mayer--Vietoris principle, the projection $p$
induces a $2$-isomorphism
\[
p^{*}:K^{*}\left( P^{*}M\right) \rightarrow K^{*}\left( S^{*}M\right)
\]
in $K$-theory of the total spaces. This readily implies that for some $N$
we have $\left[ 2^NL\right] \in \func{Im}p^{*},$ i.e. $\left[ 2^NL\right]
=\left[ 2^N\alpha ^{*}L\right] .$ This proves the first part of the lemma.

To prove the remaining part of the statement consider an arbitrary
isomorphism
\[
\sigma :2^NL\rightarrow 2^N\alpha ^{*}L.
\]
Let us extend it to the entire space
$2^N\pi ^{*}E\supset 2^NL$ in accordance with the decompositions
\begin{eqnarray*}
\widetilde{\sigma } &:&2^N\pi ^{*}E=2^NL\oplus 2^N\alpha ^{*}L\rightarrow
2^N\alpha ^{*}L\oplus 2^NL=2^N\pi ^{*}E
\end{eqnarray*}
by the formula
\begin{eqnarray*}
\widetilde{\sigma }(\xi) &=&\sigma \left( \xi \right) \oplus
\sigma \left( -\xi\right) .
\end{eqnarray*}
Then the even symbol
\[
\widetilde{\sigma }\oplus \widetilde{\sigma }^{-1}:2^{N+1}\pi
^{*}E\rightarrow 2^{N+1}\pi ^{*}E
\]
is the desired one: it has zero (topological) index
(it is even homotopic to identity), and it transforms
$2^{N+1}L$ into $2^{N+1}\alpha ^{*}L.$ The lemma is thereby proved.

\begin{lemma}
\label{factor}The quotient group $\left. \ker j\right/ \limfunc{im}i$
is a $2$-primary group.
\end{lemma}

{\em Proof\/}. Let
$$
j\left( \left[ \widehat{L}_1\right] -\left[ \widehat{L}_2\right] \right) =0,
$$
$\widehat{L}_1\subset C^\infty \left( M,E\right) ,\widehat{L}_2\subset
C^\infty \left( M,F\right)$. It follows from the definition of
the Grothendieck group $K\left({\rm Vect}(M)/\sim\right)$
that there exists an invertible even operator
\[
D:C^\infty \left( M,E\oplus E^{\prime }\right) \rightarrow C^\infty \left(
M,F\oplus E^{\prime }\right) .
\]
By considering sufficiently many copies of $E^{\prime }$, we can assume
without loss of generality that $E^{\prime }$ contains an odd
subbundle $L$ (see Lemma \ref{coke}).

The difference
\[
\left[ D\left( \widehat{L}_1\oplus \widehat{L}\right) \right] -\left[
\widehat{L}_2\oplus \widehat{L}\right]
\in K\left(\widehat{\rm Odd}/\sim\right)
\]
of odd subspaces coincides with the original element $\left[ \widehat{L}_1
\right] -\left[ \widehat{L}_2\right]$, but both subspaces $D\left(
\widehat{L}_1\oplus \widehat{L}\right)$ and $\widehat{L}_2\oplus
\widehat{L}$ lie in the same space $C^\infty \left( M,F\oplus
E^{\prime}\right)$.

In this way, the subspaces $\widehat{L}_{1,2}$
can be assumed to lie in the same space
$C^\infty\left( M,E\right)$. By Lemma \ref{fund}, for sufficiently large
$N$ there exist even invertible operators
\[
D_{1,2}:C^\infty \left( M,2^NE\right) \rightarrow C^\infty \left(
M,2^NE\right)
\]
whose principal symbols are isomorphisms of subbundles
\[
\sigma \left( D_{1,2}\right) :2^NL_{1,2}\rightarrow 2^N\alpha ^{*}L_{1,2}.
\]
Hence, the subspaces $2^N\widehat{L}_{1,2}$ and $\left( D_{1,2}\right) ^{-
1}2^N\widehat{L}_{1,2}^{\perp }$ (here $\widehat{L}^{\perp }$ denotes an
arbitrary odd complement of $\widehat{L}$) have equal principal symbols.
It follows that in the equality
\[
2^{N+1}\left[ \widehat{L}_{1,2}\right] =2^N\left[ \widehat{L}_{1,2}\oplus
\widehat{L}_{1,2}^{\perp }\right] +\left[ 2^N\widehat{L}_{1,2}\right]
-\left[ D_{1,2}^{-1}2^N\widehat{L}_{1,2}^{\perp }\right]
\in K\left(\widehat{\rm Odd}/\sim\right)
\]
the last two subspaces have equal symbols. Consequently,
by Lemma \ref{compact} we obtain
\[
\left[ 2^N\widehat{L}_{1,2}\right] -\left[ D^{-1}_{1,2}2^N
\widehat{L}_{1,2}^{\perp
}\right] \in \func{Im}i.
\]
To complete the proof of the lemma, let us show
\[
\left[ \widehat{L}_1\oplus \widehat{L}_1^{\perp }\right]=
\left[ \widehat{L}_2\oplus \widehat{L}_2^{\perp }\right] .
\]
Indeed, let $P$ and $Q$ be odd projections on $\widehat{L}_{1,2}$,
respectively. Then the even invertible operator
\begin{equation}
\begin{array}{c}
D^{\prime }:
\begin{array}{ccc}
C^\infty \left( M,E\right)  &  & C^\infty \left( M,E\right)  \\
\oplus  & \longrightarrow  & \oplus  \\
C^\infty \left( M,E\right)  &  & C^\infty \left( M,E\right)
\end{array}
, \\
D^{\prime }=\left(
\begin{array}{cc}
QP+\left( 1-Q\right) \left( 1-P\right)  & \left( 1-Q\right) P+Q\left(
1-P\right)  \\
Q\left( 1-P\right) +\left( 1-Q\right) P & QP+\left( 1-Q\right) \left(
1-P\right)
\end{array}
\right) ,
\end{array}
\label{huge1}
\end{equation}
takes $\widehat{L}_1$ to $\widehat{L}_2$. The proof of the lemma and hence
of the theorem is complete.

\begin{remark}
{\em In the proof of Lemma \ref{coke}, we constructed an odd subspace
$\widehat{L}\subset C^\infty \left(M,2^NE\right)$ in the space of sections
of the bundle $2^NE$ for an arbitrary vector bundle $E\in {\rm
Vect}\left( M\right)$. Let us show that the mapping, taking $E$ to the
subspace $\widehat{L}\oplus \widehat{L}^{\perp }$, is essentially a
splitting of the sequence (\ref{loca}). More precisely, the following
proposition holds. }
\end{remark}

\begin{proposition}
The mapping
\[
j^{\prime }:\left(
K\left({\rm Vect}(M)/\sim\right)\right)
\otimes
{\bf{}Z}\left[ \frac 12\right] \rightarrow K\left( \left. \widehat{\rm Odd}
\right/ \sim \right) \otimes {\bf{}Z}\left[ \frac 12\right]
\]
induced by the relation
\[
E\longmapsto\left[
\widehat{L}\oplus \widehat{L}^{\perp }\right] \otimes \frac 1{2^{N+1}}
\]
is a well-defined splitting of the sequence {\em (\ref{loca})}, i.e., $jj^
{\prime }=Id$.
\end{proposition}

{\em Proof\/}. One can readily show that $j^{\prime }$ is well-defined,
since
\[
\left[ \widehat{L}\oplus \widehat{L}^{\perp }\right] =\left[ \widehat{L}%
_1\oplus \widehat{L}_1^{\perp }\right]
\]
for two different decompositions of the space
$C^\infty \left( M,2^NE\right) $ into odd components, as was
already demonstrated in the end of the proof
of Lemma \ref{factor} (see~(\ref{huge1})).

The desired equality $jj^{\prime }=Id$ follows directly
from the definition. The proposition is proved.

The splitting $j^{\prime }$ of the exact sequence
(\ref{loca}) leads to the following result.

\begin{corollary}
Homomorphisms of groups
\begin{equation}
K\left( \left. \widehat{\rm Odd}\left( M\right) \right/ \sim \right) \rightarrow
{\bf{}R}  \label{star}
\end{equation}
extending the relative index ${\rm ind}:K_0\left( \left. \widehat{\rm Odd}
\right/ \sim \right) \rightarrow {\bf Z}$ of subspaces with equal
principal symbols are in a one-to-one correspondence with homomorphisms
\[
\chi : K\left({\rm Vect} (M)/\sim\right)\longrightarrow {\bf{}R;}
\]
moreover the homomorphisms {\em (\ref{star})} corresponding to integer
normalizations $\chi $ take dyadic rational values.
\end{corollary}

Let us consider the homomorphism (\ref{star}) defined by the zero
normalization $\chi \equiv 0$ in more detail. Expanding the definitions of
the Grothendieck groups, of the normalization $\chi $ and of the splitting
$j^{\prime }$ of the exact sequence (\ref{loca}), we readily obtain the
following statement.

\begin{theorem}
\label{thmn}There is a unique additive functional
\[
d:\widehat{\rm Odd}\left( M\right) \longrightarrow {\bf{}Z}\left[ \frac
12\right]
\]
with the following properties:
\begin{enumerate}
\item[{\em 1)}] invariance: for an arbitrary invertible operator $U$
with even principal symbol
$$
d\left( U\widehat{L}\right) =d\left( \widehat{L}\right);
$$ 
\item[{\em 2)}] relative dimension:
for two subspaces $\widehat{L}_1$ and $\widehat{L}_2$ with equal principal
symbols
$$
{\rm ind}\left(
\widehat{L}_1,\widehat{L}_2\right) =d\left( \widehat{L}_1\right) -d\left(
\widehat{L}_2\right);
$$
\item[{\em 3)}] complement:
for an arbitrary odd subspace $\widehat{L}$
$$
d\left( \widehat{L}\right) +d\left( \widehat{L}^{\perp}\right) =0.
$$
\end{enumerate}
\end{theorem}

\begin{remark}
The first property of the functional $d$ implies its homotopy invariance.
\end{remark}

Indeed, it was shown in the proof of lemma \ref{compact} that for a smooth
homotopy of odd projections $P_t$ corresponding to a family of subspaces
$\widehat{L}_t$ it is possible to construct a family of invertible
operators $U_t$ such that $\widehat{L}_t=U_t\widehat{L}_0$. Thus, by the
invariance property of the functional $d$, we obtain
\[
d\left( \widehat{L}_t\right) ={\rm const}.
\]

\begin{remark}\label{impor}
{\em The problem of finding an explicit formula for $d$ is not simple. The
previous theorem only shows that this invariant can be expressed in terms
of the relative index of subspaces and the invertible even operator $D$
constructed by the principal symbol from Lemma \ref{fund},
\[
D:C^\infty \left( M,2^NE\right) \rightarrow C^\infty \left( M,2^NE\right)
,\quad \sigma \left( D\right) :2^NL\rightarrow 2^N\alpha ^{*}L,
\]
by the formula
\[
d\left( \widehat{L}\right) =\frac 1{2^{N+1}}{\rm ind}
\left( D2^N\widehat{L},2^N%
\widehat{L}^{\perp }\right) .
\]
In section \ref{undeta} we obtain a formula for $d$ in terms of the $\eta
$-invariant of Atiyah--Patodi--Singer. }
\end{remark}

\section{The index of operators in subspaces}

Let us recall the notion of a pseudodifferential operator acting in
subspaces \cite{SaSt1}.

A pseudodifferential operator
\[
D:C^\infty \left( M,E\right) \longrightarrow C^\infty \left( M,F\right)
\]
of order $m$ is said to act in subspaces
\[
\widehat{L}_1\subset C^\infty \left( M,E\right) ,\quad \widehat{L}_2\subset
C^\infty \left( M,F\right)
\]
if
\[
D\widehat{L}_1\subset \widehat{L}_2\text{.}
\]
For the restriction
\begin{equation}
D:\widehat{L}_1\longrightarrow \widehat{L}_2  \label{subsp}
\end{equation}
of the operator $D$ to the subspaces, there is a criterion for the
Fredholm property of the closure of this operator with respect to Sobolev
norms (see \cite{ScSS18}).

\begin{proposition}
An operator
\[
D:\widehat{L}_1\longrightarrow \widehat{L}_2,
\]
where
\[
 \widehat{L}_1\subset
H^s\left( M,E\right) ,\quad \widehat{L}_2\subset H^{s-m}\left( M,F\right),
\]
is Fredholm if and only if it is elliptic, i.e., its principal symbol
\[
\sigma \left( D\right) :L_1\longrightarrow L_2
\]
is an isomorphism of the vector bundles $\ L_1$ and $L_2.$
\end{proposition}

For the case of operators acting in odd subspaces on even-dimensional
manifolds, which were studied in the previous section,
let us obtain an index formula it terms of the functional $d$.

\begin{theorem}
For an operator
\[
D:\widehat{L}_1\longrightarrow \widehat{L}_2,
\]
where
\[
\widehat{L}_{1,2}\in
\widehat{\rm Odd}\left( M^{ev}\right) ,\quad \widehat{L}_1\subset C^\infty
\left( M,E\right) ,\widehat{L}_2\subset C^\infty \left( M,F\right),
\]
the following index formula is valid:
\begin{equation}
{\rm ind}\left( D,\widehat{L}_1,\widehat{L}_2\right) =\frac 12{\rm ind}\widetilde{D}%
+d\left( \widehat{L}_1\right) -d\left( \widehat{L}_2\right) ,  \label{oba}
\end{equation}
where
$$
\widetilde{D}: C^\infty
\left( M,E\right)\longrightarrow
C^\infty\left( M,F\right)
$$
is an elliptic operator with principal symbol
\[
\begin{array}{cccc}
\sigma \left( \widetilde{D}\right) : & \pi ^{*}E & \rightarrow  & \pi ^{*}F,
\\
& L_1\oplus \alpha ^{*}L_1 & \stackrel{\sigma \left( D\right) \oplus \alpha
^{*}\sigma \left( D\right) }{\longrightarrow } & L_2\oplus \alpha ^{*}L_2.
\end{array}
\]
\end{theorem}

{\em Proof\/}. Modifying, if necessary, the subspaces $\widehat{L}_{1,2}$
by finite-dimensional subspaces, we can assume that $D$ is an isomorphism.
\footnote{One can obtain this, for example, by deleting the kernel of the
operator $D$ from the subspace $\widehat{L}_{1}$ and the cokernel of this
operator from $\widehat{L}_{2}$.} Let us choose the operator
$\widetilde{D}$ in such a way that it coincides with the original operator
$D$ on the subspace $\widehat{L}_{1}$. To be definite, we can also assume
that ${\rm ind}\widetilde{D}\geq 0$ and the following relations are valid:
\begin{equation}
\limfunc{coker}\widetilde{D}=0,\quad
\ker \widetilde{D}\subset \widehat{L}_1^{\perp }.
\label{treb}
\end{equation}
(In the case ${\rm ind}\widetilde{D}< 0$, the index formula can be proved
for the inverse of $D$.) Let us verify the index formula (\ref{oba}) in
this case.

On the one hand, by virtue of conditions (\ref{treb}), we have
\begin{equation}
{\rm ind}\left( D,\widehat{L}_1,\widehat{L}_2\right) =0,\;{\rm ind}\widetilde{D}=\dim
\ker \widetilde{D}.  \label{uzhe}
\end{equation}
Let us now calculate the difference $d\left( \widehat{L}_1\right)
-d\left( \widehat{L}_2\right)$. To this end, we consider
an invertible even operator $U$
\[
U:\begin{array}{ccc}
C^\infty \left( M,E\right) &  & C^\infty \left( M,E\right) \\
\oplus & \longrightarrow & \oplus \\
C^\infty \left( M,F\right) &  & C^\infty \left( M,F\right)
\end{array}
,\;U=\left(
\begin{array}{cc}
\vspace{0.2cm}
P_{\ker \widetilde{D}} & \widetilde{D}^{-1} \\
\widetilde{D} & 0
\end{array}
\right) ,
\]
where $P_{\ker \widetilde{D}}$ is the projection on the
(finite-dimensional) kernel of $\widetilde{D}$  and
the operator $\widetilde{D}^{-1}$ maps into the orthogonal
complement to this kernel.

It is easy to check that the operator $U$ isomorphically maps
\begin{equation}
\widehat{L}_1\oplus \widetilde{D}
\widehat{L}_1^{\bot }\quad \text{onto\quad }%
\left( \widehat{L}_1^{\bot }\ominus \ker \widetilde{D}\right) \oplus
\widehat{L}_2.
\label{laba}
\end{equation}
It follows from (\ref{treb}) that $\widetilde{D}%
\widehat{L}_1^{\bot }= \widehat{L}_2^{\bot }$. By applying
the functional $d$ to the subspaces in (\ref{laba}), we obtain
\[
d\left( \widehat{L}_1\right) -d\left( D\widehat{L}_1\right) =-d\left(
\widehat{L}_1\right) -\dim \ker \widetilde{D}+d\left( \widehat{L}_2\right) .
\]
Hence
\[
d\left( \widehat{L}_1\right) -d\left( \widehat{L}_2\right) =-\frac 12\dim
\ker \widetilde{D}=-\frac 12{\rm ind}\widetilde{D},
\]
which, together with relations (\ref{uzhe}), proves the desired formula
\[
{\rm ind}\left( D,\widehat{L}_1,\widehat{L}_2\right) =
\frac 12{\rm ind}\widetilde{D}+d\left( \widehat{L}_1\right)
 -d\left( \widehat{L}_2\right) .
\]
The proof is complete.

\begin{corollary}
The following index formula is valid for an operator acting between a
subspace and a space:
\begin{equation}
{\rm ind}\left( D,\widehat{L},C^\infty \left( M,F\right) \right) =\frac 12{\rm ind}%
\widetilde{D}+d\left( \widehat{L}\right) ,  \label{gladk}
\end{equation}
where the principal symbol of the operator $\widetilde{D}$
\[
\widetilde{D}:C^\infty \left( M,E\right) \rightarrow C^\infty \left(
M,F\oplus F\right)
\]
is given by the formula
\[
\sigma \left( \widetilde{D}\right) =\sigma \left( D\right) \oplus \alpha
^{*}\sigma \left( D\right) .
\]
\end{corollary}

{\em Proof\/}. By Lemma \ref{coke}, for sufficiently
large $N$ the space $C^\infty \left( M,2^NE\right) $ can be
decomposed into a direct sum of two odd subspaces. Formula
(\ref{gladk}) then readily follows from (\ref{oba}).

\section{Dimension and the $\eta $-invariant\label{undeta}}

Recall \cite{Gil7} that a classical pseudodifferential
operator $A$ of a positive integer order $m$ is said to be
{\em admissible\/} if its complete symbol
\[
a\left( x,\xi \right) \sim a_m\left( x,\xi \right) +a_{m-1}\left( x,\xi
\right) +...
\]
satisfies the parity conditions
\[
a_\alpha \left( x,-\xi \right) =\left( -1\right) ^\alpha a_\alpha \left(
x,\xi \right) ,\qquad \forall \xi \neq 0,\quad x,\quad\alpha=m,m-1,m-2,\ldots
\]
(the admissibility of the operator is independent of the choice of
the coordinate system where the complete symbol is considered).

\begin{theorem}
For an admissible self-adjoint elliptic operator $A$ of odd order with
nonnegative spectral subspace $\widehat{L}$ on an even-dimensional
manifold, the following formula is valid:
\[
 d\left( \widehat{L}\right)=\eta \left( A\right).
\]
\end{theorem}

{\em Proof\/}. The fact that the operator $A$ is admissible
implies that its principal symbol is odd,
\[
\alpha ^{*}\sigma \left( A\right) =-\sigma \left( A\right) .
\]
Hence, it follows from Proposition \ref{opisub} that the nonnegative
spectral subspace $\widehat{L}$ is odd (and even orthogonally odd). By
Lemma \ref{fund}, there exists an even isomorphism $\sigma $
\[
\sigma :2^N\pi ^{*}E\rightarrow 2^N\pi ^{*}E
\]
with zero (topological) index that takes the subbundle $2^NL$ to
$2^N\alpha ^{*}L$. Consider an invertible even-order operator $U_0$ with
principal symbol $\sigma$. We can also assume that the operator $U_0$ is
admissible. An operator with these properties can be constructed in the
following way. By the usual method, we construct an admissible
pseudodifferential operator with principal symbol $\sigma$, using local
charts on the manifold. Then, modifying the operator thus obtained by a
finite rank operator, we construct an invertible operator. Let us
construct a unitary operator $U$ by the formula
\[
U=U_0\left( \sqrt{U_0^{*}U_0}\right) ^{-1}.
\]
It isomorphically translates the subspace $2^N\widehat{L}$ to a subspace
with the principal symbol equal to $2^N\alpha ^{*}L.$ Thus, we obtain the
following expression for the functional $d$ via the relative index of
subspaces (see Remark \ref{impor})
\begin{equation}
2^{N+1}d\left( \widehat{L}\right) ={\rm ind}\left( U2^N\widehat{L},2^N\widehat{L}%
^{\bot }\right) .  \label{ford}
\end{equation}
On the other hand, the $\eta$-invariant of an elliptic self-adjoint
operator is determined by the spectrum of the operator.
Therefore, the unitary property of the operator $U$ implies
\[
\eta \left( 2^NA\right) =\eta \left( U2^NAU^{-1}\right) .
\]
The operators $-2^NA$ and $U2^NAU^{-1}$ have equal principal symbols and
are admissible operators of odd order on an even-dimensional manifold.
Hence, the difference of their $\eta $-invariants is equal to the relative
index of the nonnegative spectral subspaces for these self-adjoint
operators:
\[
\eta \left( -2^NA\right) -\eta \left( U2^NAU^{-1}\right) ={\rm ind}\left( \widehat{%
L}_{+}\left( -2^NA\right) ,\widehat{L}_{+}\left( U2^NAU^{-1}\right) \right) .
\]
Thus, we express the $\eta$-invariant in terms of the relative index by
the formula
\[
2^{N+1}\eta \left( A\right) ={\rm ind}\left( U2^N\widehat{L},2^N\widehat{L}^{\bot
}\right) .
\]
Comparing the last formula with the expression
(\ref{ford}), we obtain the desired identity
\[
d\left( \widehat{L}\right)=\eta \left( A\right).
\]

\begin{remark}
\label{keta}
{\em
An odd subspace $\widehat{L}\subset C^\infty \left( M,E\right) $ can be
realized as the nonnegative spectral subspace of an admissible elliptic
self-adjoint operator of odd order if and only if it is orthogonally odd,
the orthogonal projection $P$ onto it is a classical pseudodifferential
operator with symbol
\[
\sigma \left( P\right) \left( x,\xi \right) \sim \sum_{j\leq 0}p_j\left(
x,\xi \right) ,
\]
and the homogeneous lower-order terms of the complete symbol
have the following symmetry for $j<0$
\begin{equation}
p_j\left( x,-\xi \right) =\left( -1\right) ^{j+1}p_j\left( x,\xi \right) .
\label{alagil}
\end{equation}
}
\end{remark}
Indeed, for an invertible self-adjoint elliptic operator $A$,
the positive spectral projection is equal to
\[
P_{+}\left( A\right) =\frac{1+A\left| A\right| ^{-1}}2.
\]
From this equality, using the calculus of complete symbols
of pseudodifferential operators, we can verify the necessity
of the condition (\ref{alagil}).

Conversely, for a classical pseudodifferential orthogonal
projection $P$ satisfying (\ref{alagil}), the desired
self-adjoint operator can be defined by the formula
\[
A=\sqrt[4]{\Delta }\left( 2P-1\right) \sqrt[4]{\Delta },
\]
where $\Delta $ is a positive Laplacian in the bundle $E$.

\section{The index of boundary value problems with odd projections}

Boundary value problems for general elliptic operators, that is, operators
that in general violate the Atiyah-Bott condition \cite{AtBo2}, were
introduced in \cite{ScSS18}. These boundary value problems have the form
\begin{equation}
\left\{
\begin{array}{l}
Du=f\in H^{s-m}\left( M,F\right),\quad u\in H^s\left( M,E\right), \\
Bu=g\in \widehat{L}\subset H^\delta \left( \partial M,G\right) ,
\end{array}
\right.   \label{bvp1}
\end{equation}
where the subspace $\widehat{L}$ is the range of a pseudodifferential
projection and $B$ is some boundary operator.

In the paper \cite{SaSt1} we obtained a reduction of a general boundary
value problem of the form (\ref{bvp1}) to the so-called {\em spectral
boundary value problem\/} \cite{NScSS3} for a first-order operator of a
special form in a neighborhood of the boundary of manifold $M$. Thus, in
the present paper we can restrict ourselves to the index problem for
spectral boundary value problems only. Let us describe these problems in
more detail.

On a smooth compact manifold $M$ with boundary
$X=\partial M$, consider an elliptic operator
\[
D:C^\infty \left( M,E\right) \longrightarrow C^\infty \left( M,F\right)
\]
that is a first-order elliptic pseudodifferential operator
outside some neighborhood of the boundary, while in the collar
neighborhood
\[
U_{\partial M}=\left[ 0,1\right) \times \partial M
\]
with the normal coordinate $t$ to the boundary, it has the form
\begin{equation}
\left. D\right| _{U_{\partial M}}=\chi \left( t\right) \left[ -i\frac%
\partial {\partial t}-iA\left( t\right) \right] ,  \label{aleph}
\end{equation}
for a smooth isomorphism $\chi \left( t\right) $ of the
vector bundles $E$ and $F$, where $A\left( t\right) $ is a
smooth family of first-order pseudodifferential operators
on the boundary with principal symbols
$\sigma \left( A\left(t\right) \right) \left( \xi ^{\prime }\right) $
that have no pure imaginary eigenvalues for
$\left| \xi^{\prime }\right| \neq 0$. Moreover, it can be assumed
that the principal symbol of the operator $A\left( 0\right) $
on the unit cotangent spheres to the boundary
$S^{*}X$ has only two eigenvalues $\pm 1$,
\begin{equation}
\sigma \left( A\left( 0\right) \right) ^2=1\qquad \text{for }%
\left| \xi ^{\prime }\right| =1,  \label{mu}
\end{equation}
and has a basis of eigenvectors.

\begin{remark}
{\em
The operator $D$ in  (\ref{aleph}) is in  the general case an elliptic
operator with {\em continuous} symbol in the sense of \cite{ReSc1}. }
\end{remark}

\begin{definition}
A spectral boundary value problem
{\em for an elliptic operator $D$ on a manifold $M$ with
boundary is a system of equations of the form
\begin{equation}
\left\{
\begin{array}{ll}
Du=f & \in H^{s-1}\left( M,E\right), u\in H^{s}(M,F),  \\
P\left. u\right| _{\partial M}=g & \in \widehat{L}=\func{Im}P\subset
H^{s-1/2}\left( \partial M,\left. E\right| _{\partial M}\right)
\end{array}
\right.   \label{spbvp}
\end{equation}
(see \cite{SaSt1}), where an admissible subspace $\widehat{L}$ on the
boundary $\partial M$ has the symbol equal to the nonnegative spectral
subspace of the symbol $\sigma \left( A\left( 0\right) \right) $. }
\end{definition}

Let us denote  boundary value problems of the form
(\ref{spbvp}) by $\left( D,\widehat{L}\right) .$

Consider a spectral boundary value problem (\ref{spbvp}) such that
the subspace $\widehat{L}$ is odd. In this case, the principal symbol
of $D$ can be extended to the double
\[
2M\stackrel{{\rm def}}{=}M\bigcup\limits_{\partial M}M
\]
of the manifold $M$ by the formula
\[
\sigma \left( D\right) \bigcup\limits_{\left. S^{*}M\right| _{\partial
M}}\alpha ^{*}\sigma ^{-1}\left( D\right).
\]
Let us verify the continuity of the symbol obtained at the place of gluing
of the manifolds, i.e. let us show the coincidence of the symbols
$\sigma \left( D\right) $ and $\alpha ^{*}\sigma^{-1} \left( D\right)$
on the bundle ${\left. S^{*}M\right| _{\partial M}}$ of cotangent
spheres to the manifold $M$ over the boundary.

Let us denote the dual variables
to the normal and tangent variables on the boundary
$\partial M$ by $\tau $ and $\xi ^{\prime }$, respectively.
The principal symbol of the operator
$A\left( 0\right) $ is denoted by $\sigma \left( \xi ^{\prime }\right)$.
Formula (\ref{aleph}) then implies
(for the isomorphism $\chi(0)=$Id for brevity)
\[
\sigma \left( D\right) \left( \tau ,\xi ^{\prime }\right) =\tau -i\sigma
\left( \xi ^{\prime }\right) .
\]
On the other hand, since $\sigma \left( \xi ^{\prime }\right) $ is odd, we
have
\[
\sigma \left( D\right) \left( \tau ,-\xi ^{\prime }\right) =\tau +i\sigma
\left( \xi ^{\prime }\right).
\]
Hence, by (\ref{mu}) we obtain
\[
\sigma ^{-1}\left( D\right) \left( \tau ,-\xi ^{\prime }\right) =\tau
-i\sigma \left( \xi ^{\prime }\right)
\]
on the unit cospheres ${\left. S^{*}M\right| _{\partial
M}}:\tau^{2}+|\xi'|^{2}=1$. In this way, we establish the equality
\[
\sigma \left( D\right) \left( \tau ,\xi ^{\prime }\right) =\sigma
^{-1}\left( D\right) \left( \tau ,-\xi ^{\prime }\right).
\]
This equation is readily seen to be equivalent to the condition of the
continuous gluing for the elliptic symbols $\sigma \left( D\right) $ and
$\alpha ^{*}\sigma ^{-1}\left( D\right) $ that come from the two halves
of the double $2M$.

An elliptic operator constructed from the symbol thus obtained is
denoted by $\widetilde{D}$:
\begin{equation}
\sigma \left( \widetilde{D}\right) =\sigma \left( D\right)
\bigcup\limits_{\left. S^{*}M\right| _{\partial M}}\alpha ^{*}\sigma
^{-1}\left( D\right) .  \label{omega}
\end{equation}

On an even-dimensional manifold $M$, the procedure just described leads to
one of the terms in the index formula for boundary value problems with odd
subspaces. We state this formula in the form of a theorem.

\begin{theorem}
For a spectral boundary value problem $\left( D,\widehat{L}\right) $ with
an odd subspace $\widehat{L}\in \widehat{\limfunc{Odd}}
\left( M^{ev}\right)$, the following index formula is valid
\begin{equation}
{\rm ind}\left( D,\widehat{L}\right) =\frac 12\,{\rm ind}\,
\widetilde{D}-d\left( \widehat{L}\right) ,
\label{thind}
\end{equation}
where $\widetilde{D}$ is an elliptic operator on the double $2M$ with
principal symbol {\em (\ref{omega})} and the functional $d$ is defined in
Section {\em \ref{dime}}.
\end{theorem}

{\em Proof\/}. By analogy with the case of even boundary value problems
\cite{SaSt1}, the basic idea of the proof is to make a reduction to a
boundary value problem that is classical in a certain sense.

We check, first of all, that the index formula in question is true for the
case in which the subspace $\widehat{L}$ is a sum of two complementary
subspaces,
\begin{equation}
\widehat{L}=\widehat{L}_0\oplus \widehat{L}_0^{\perp },\qquad \widehat{L}%
_0\subset C^\infty \left( \partial M,G\right) .  \label{Omeg}
\end{equation}
By virtue of the complement property (see Theorem \ref{thmn}), we obtain
$d\left( \widehat{L}\right) =0$. Thus, the check is reduced to
finding
out whether
\begin{equation}
{\rm ind}\left( D,\widehat{L}_0\oplus \widehat{L%
}_0^{\perp }\right)=\frac 12\limfunc{ind}\widetilde{D}.
\label{star2}
\end{equation}
To this end, we take a closer look at the spectral boundary value problem
$\left( D,\widehat{L}_0\oplus \widehat{L}_0^{\perp }\right)$. This
problem, by virtue of the isomorphism
$$
\begin{array}{c}
B_0:\widehat{L}_0\oplus \widehat{L}_0^{\perp }\longrightarrow C^\infty
\left( \partial M,G\right) , \\
u\oplus v \longrightarrow u+v,
\end{array}
$$
can be viewed as a classical boundary value problem. The cohomological
index formula for classical boundary value problems then implies
(cf. the proof of (\ref{A3}))
$$
{\rm ind}_t\left( \sigma \left( D\right) ,L_0\oplus \alpha ^{*}L_0\right)
=-{\rm ind}_t\left( \alpha ^{*}\left[ \sigma \left( D\right)
,L_0\oplus \alpha^{*}L_0\right] \right).
$$
Hence,
\[
2\cdot{\rm ind}\left( D,\widehat{L}_0\oplus \widehat{L}_0^{\perp }\right)
={\rm ind}_t\left(\sigma \left( D\right) ,L_0\oplus \alpha ^{*}L_0  \right)
 -{\rm ind}_t\left(\alpha ^{*}\left[ \sigma \left( D\right)
,L_0\oplus \alpha^{*}L_0\right] \right).
\]
Finally, the additivity of the topological index gives
\[
\!{\rm ind}_t\left(\sigma \left( D\right) ,L_0\oplus \alpha ^{*}L_0  \right)
\! -\!{\rm ind}_t\left(\alpha ^{*}\left[ \sigma \left( D\right)
,L_0\oplus \alpha^{*}L_0\right] \right)
={\rm ind}_t\left( \sigma \left( D\right) \cup \alpha%
^{*}\sigma ^{-1}\left( D\right) \right) ={\rm ind}\widetilde{D}.
\]
Thus, we have proved the index formula for the decomposition (\ref{Omeg}).
Let us now verify that the index formula is also valid for the case in
which the decomposition (\ref{Omeg}) takes place for the symbols of
subspaces,
\begin{equation}
\label{xdeco}
L=L_0\oplus \alpha ^{*}L_0.
\end{equation}
Indeed, by the logarithmic property of the index, we obtain
\begin{equation}
{\rm ind}\left( D,\widehat{L}\right) ={\rm ind}\left( D,\widehat{L}_0\oplus \widehat{L}%
_0^{\perp }\right) +{\rm ind}\left( \widehat{L}_0\oplus \widehat{L}_0^{\perp },%
\widehat{L}\right).  \label{verif}
\end{equation}
On the one hand,
\[
{\rm ind}\left( \widehat{L}_0\oplus \widehat{L}_0^{\perp },\widehat{L}\right)
=d\left( \widehat{L}_0\oplus \widehat{L}_0^{\perp }\right) -d\left( \widehat{%
L}\right) =-d\left( \widehat{L}\right) ,
\]
and on the other hand, we have already proved that
\[
{\rm ind}\left( D,\widehat{L}_0\oplus \widehat{L}_0^{\perp }\right) =\frac 12{\rm ind}%
\widetilde{D}.
\]
Substituting these two results into (\ref{verif}), we obtain
the index formula in this case as well.

To complete the proof of the index formula (\ref{thind}),
let us reduce a general boundary value problem
$\left( D,\widehat{L}\right) $ to the form already considered.
For this purpose, let us apply the invertible even operator
\[
U:C^\infty \left( \partial M,2^NE\right) \longrightarrow C^\infty \left(
\partial M,2^NE\right) ,
\]
given by Lemma \ref{fund} of Section $2$,
with principal symbol $u$ that establishes an isomorphism of bundles
$$
u: 2^NL\longrightarrow\alpha ^{*}2^NL.
$$
This operator enables us to deform the odd subspace $2^{N+2}\widehat{L}$
to a subspace with principal symbol having a decomposition of the
form (\ref{xdeco})
\[
2^{N+1}\left( L\oplus \alpha ^{*}L\right),
\]
by the formula
\begin{equation}
\begin{array}{c}
U_\varphi :
\begin{array}{ccc}
C^\infty \left( \partial M,2^{N+2}E\right)  & \longrightarrow  & C^\infty
\left( \partial M,2^{N+2}E\right)
\end{array}
,\qquad \varphi \in \left[ 0,\pi /2\right] , \\
U_\varphi =1_{2^{N+2}E}\oplus \left(
\begin{array}{cc}
\cos \varphi  & -U^{-1}\sin \varphi  \\ %
U\sin \varphi  & \cos \varphi %
\end{array}
\right) .
\end{array}
\label{huge2}
\end{equation}
In \cite{SaSSc1} it is shown that the homotopy of subspaces occuring in
the boundary value problem can be lifted to a homotopy of spectral
boundary value problems by deforming the operator $D$ in a neighborhood of
the boundary. Hence, the spectral boundary value problem
\[
\left( 2^{N+2}D,2^{N+2}\widehat{L}\right)
\]
is homotopic to a spectral boundary value problem with principal symbol
of the form (\ref{star2})
\[
\left( \sigma \left( D^{\prime }\right) ,2^{N+1}\left( L\oplus \alpha
^{*}L\right) \right) .
\]
This completes the proof of the theorem, since the components of the index
formula (\ref{thind}) are homotopy invariant.

\section{Example}
On a compact smooth manifold $M$ without boundary consider an elliptic
pseudodifferential operator $D$
$$
D:C^\infty(M,E)\longrightarrow C^\infty(M,F).
$$
From the operator $D$ we construct an elliptic
self-adjoint operator $A$ by the formula
$$
A=\left(\begin{array}{cc} 0 & D^* \\ D & 0 \end{array} \right):
\begin{array}{c} C^\infty(M,E) \\ \oplus \\ C^\infty(M,F) \end{array}
\longrightarrow
\begin{array}{c} C^\infty(M,E) \\ \oplus \\ C^\infty(M,F) \end{array}.
$$
Let $\widehat{L}_+(A)$ be the subspace generated by the
eigenvectors of $A$ corresponding to nonnegative eigenvalues;
\begin{equation}
\label{dek1}
\widehat{L}_+(A)\subset C^\infty(M,E)\oplus C^\infty(M,F).
\end{equation}
The orthogonal projection on the second term in this formula
will be denoted by $P_0$.
\begin{proposition}
\label{dirind}
The operator in subspaces
$$
P_0: \widehat{L}_+(A)\longrightarrow  C^\infty(M,F)
$$
is elliptic, and its index can be calculated by the formula
$$
\ker P_0=\ker D\oplus 0, {\rm coker}P_0=0,\quad
{\rm ind}(P_0,\widehat{L}_+(A),C^\infty(M,F))={\rm dim}\ker D.
$$
\end{proposition}

{\em Proof\/}. First, let us calculate the subspace $\widehat{L}_+(A)$ in
terms of the decomposition (\ref{dek1}). A simple calculation shows that
this subspace is generated by the following eigenvectors of the operator $A$:
\begin{equation}
\label{eigensub}
\left(
\begin{array}{c}
u \\
\frac{Du}{\sqrt{\lambda}}
\end{array}
\right),
\left(\begin{array}{c} u' \\ 0 \end{array} \right),\quad u'\in{\rm ker} D,
\left(\begin{array}{c} 0 \\ v \end{array} \right),\quad v\in{\rm ker} D^*,
\end{equation}
where $u$ is an eigenvector of the nonnegative operator
$D^*D$ corresponding to a positive eigenvalue $\lambda$.
The statement we need to prove follows directly from this formula.

Now we calculate the index of the operator $P_0$, using the
index formula for operators in subspaces (see formula (\ref{gladk})).
To this end, we suppose that the operator $D$ is,
for example, an odd-order differential operator
on an even-dimensional manifold $M$.

Let us calculate the terms in the index formula
\begin{equation}
{\rm ind}\left( P_0,\widehat{L}_+(A),C^\infty \left( M,F\right) \right)
=\frac 12{\rm ind}%
\widetilde{P}_0+\eta\left(A\right).
\end{equation}
The elliptic operator $\widetilde{P}_0$
$$
\widetilde{P}_0:C^\infty \left( M,E\oplus F\right)\longrightarrow
C^\infty \left( M,F\oplus F\right),
$$
by definition, has the principal symbol
$$
\sigma(\widetilde{P}_0)(\xi)=
\sigma({P}_0)(\xi)\oplus\sigma({P}_0)(-\xi):
L_+(A)\oplus L_-(A)\longrightarrow\pi^*(F\oplus F).
$$
This, in compliance with formula (\ref{eigensub}), permits us to
rewrite this symbol in the form
$$
\sigma(\widetilde{P}_0)=
\frac{\sigma(D)}{\sqrt{\sigma(D^*)\sigma(D)}}\oplus 1:
\pi^*(E\oplus F)\longrightarrow\pi^*(F\oplus F).
$$
From this we obtain
$$
{\rm ind}\widetilde{P}_0={\rm ind}D.
$$

We finally calculate the $\eta$-invariant of the operator $A$. This
operator has a spectral symmetry (i.e., its eigenvalues lie symmetrically
on the real line with respect to zero, and eigenvalues of opposite signs
have equal multiplicities). Thus, the $\eta$-function of the operator is
constant:
$$
\eta(A)(s)=\frac12\left({\rm dim\;ker} A+
\sum_{\lambda_i\ne 0}|\lambda_i|^{-s}{\rm sgn}\lambda_i
\right)
\equiv
\frac{{\rm dim\;ker} A}2.
$$
For the $\eta$-invariant we obtain
$$
\eta(A)=\eta(A)(0)=\frac{{\rm dim\;ker} D+{\rm dim\;ker} D^*}2=
-\frac{{\rm ind} D}2+{\rm dim\;ker} D.
$$
Hence, the right-hand side of the index formula is equal to
$$
\frac 12{\rm ind}\widetilde{P}_0+\eta\left(A\right)=
\frac{{\rm ind} D}2-\frac{{\rm ind} D}2+{\rm dim\;ker} D=
{\rm dim\;ker} D,
$$
which coincides with the index computation for the operator $P_0$ from
Proposition \ref{dirind}.

\appendix

\section{Appendix}

Consider two subspaces $\widehat{L}_{1,2}$ with equal principal symbols.
The {\em relative index of subspaces\/} $\widehat{L}_1$ and $\widehat{L}_2
$ is, by definition, the index of the following Fredholm operator defined
by the projections on these subspaces:
\[
{\rm ind}\left( \widehat{L}_1,\widehat{L}_2\right) =
{\rm ind}\left( P_2:\func{Im}%
P_1\rightarrow \func{Im}P_2\right) ,\quad \text{for }\widehat{L}%
_{1,2}=\func{Im}P_{1,2}.
\]
The relative index is independent of the choice of projections
$P_{1,2}$ on the corresponding subspaces. It has the following
geometric interpretation.

Consider the Grassmanian ${\rm Gr}(L)$ of all admissible subspaces with
symbol $L$ (we assume that the symbol $L$ is nontrivial, i.e., $0<\dim L
<\dim E$). The following theorem was proved in \cite{Wojc1}.

\begin{theorem}
\label{thwojc}
The Grassmanian ${\rm Gr}(L)$ is a classifying space for $K$-theory.
\end{theorem}

On the other hand, it is known by Atiyah--J{\"a}nich theorem \cite{Ati2}
that the classifying space for $K$-theory is the space of Fredholm
operators in a given separable Hilbert space $H$. That is, for a space
$X$ there is an isomorphism of groups
\begin{equation}
{\rm ind}:\left[X,{\rm Fred}(H)\right]\longrightarrow K(X),
\end{equation}
where the map is given by the index of a family of Fredholm operators.

Both realizations of $K$-theory (i.e., by means of Theorem \ref{thwojc} or
by the Atiyah--J{\"a}nich theorem) are linked by a simple formula namely,
by the relative index of subspaces. Specifically, by fixing a certain
subspace $\widehat{L}_0$ with principal symbol $L$, we assign a Fredholm
operator to any subspace $\widehat{L}\in{\rm Gr}(L)$ by the formula
$$
P_0:\widehat{L}{\longrightarrow}\widehat{L}_0,\quad
\widehat{L}_0={\rm Im}P_0.
$$
Thus, we define a mapping to the set of Fredholm operators in Hilbert spaces
$$
{\rm Gr}(L)\longrightarrow{\rm Fred}.
$$
We point out that in this case we obtain a Fredholm operator that,
depending on the point of the Grassmanian $\widehat{L}\in{\rm Gr}(L)$,
acts in different Hilbert spaces.\footnote{Meanwhile, by the theorem of
Kuiper \cite{Kui1} on the contractibility of the unitary group of the
Hilbert space, the resulting operators can be reduced to operators in a
fixed Hilbert space.}
In terms of the Grassmanian ${\rm Gr}(L)$, the isomorphism map
$$
\left[X,{\rm Gr}(L)\right]\stackrel{\beta}{\longrightarrow} K(X)
$$
and its inverse are given by the formulas
$$
\beta\left(\widehat{L}(x)\right)={\rm ind}\left(
P_0:\widehat{L}(x)
{\longrightarrow}\widehat{L}_0
\right)\in K(X),
$$
where $\widehat{L}(x)$ denotes a family of subspaces with symbol
equal to $L$ and with parameter $x\in X$. The inverse mapping is
given by the formula
$$
\beta^{-1}\left([E]-[F]\right)=\left[(\widehat{L}_0\ominus F)\oplus E\right]
\in [X, {\rm Gr}(L)],
$$
where the family of subspaces denoted by
$(\widehat{L}_0\ominus F)\oplus E$ is constructed as follows:
since $\widehat{L}_0$ is finite-dimensional and finite-codimensional,
we embed the vector bundle $F$ in the trivial vector
bundle over $X$ with fiber $\widehat{L}_0$, while $E$ is embedded
in the bundle with fiber $\widehat{L}^\perp_0$, which is
an orthogonal complement to $\widehat{L}_0$.



\begin{flushright}
\em
Moscow  State University
\end{flushright}
\end{document}